\documentclass[10pt,a4paper]{article}
\usepackage{amssymb}
\usepackage{amscd}
\usepackage{latexsym}
\usepackage{amsmath}
\usepackage{amsfonts}
\usepackage{amscd}
\usepackage{xcolor}
\usepackage{ulem}
\begin{document}

\newtheorem{theorem}{Theorem}[section]
\newtheorem{remark}[theorem]{Remark}
\newtheorem{mtheorem}[theorem]{Main Theorem}
\newtheorem{bbtheo}[theorem]{The Strong Black Box}
\newtheorem{observation}[theorem]{Observation}
\newtheorem{proposition}[theorem]{Proposition}
\newtheorem{lemma}[theorem]{Lemma}
\newtheorem{testlemma}[theorem]{Test Lemma}
\newtheorem{mlemma}[theorem]{Main Lemma}
\newtheorem{note}[theorem]{{\bf Note}}
\newtheorem{steplemma}[theorem]{Step Lemma}
\newtheorem{corollary}[theorem]{Corollary}
\newtheorem{notation}[theorem]{Notation}
\newtheorem{example}[theorem]{Example}
\newtheorem{definition}[theorem]{Definition}

\renewcommand{\labelenumi}{(\roman{enumi})}
\def\Pf{\smallskip\goodbreak{\sl Proof. }}

\def\Fin{\mathop{\rm Fin}\nolimits}
\def\br{\mathop{\rm br}\nolimits}
\def\fin{\mathop{\rm fin}\nolimits}
\def\Ann{\mathop{\rm Ann}\nolimits}
\def\Aut{\mathop{\rm Aut}\nolimits}
\def\End{\mathop{\rm End}\nolimits}
\def\Tor{\mathop{\rm Tor}\nolimits}
\def\bfb{\mathop{\rm\bf b}\nolimits}
\def\bfi{\mathop{\rm\bf i}\nolimits}
\def\bfj{\mathop{\rm\bf j}\nolimits}
\def\df{{\rm df}}
\def\bfk{\mathop{\rm\bf k}\nolimits}
\def\bEnd{\mathop{\rm\bf End}\nolimits}
\def\iso{\mathop{\rm Iso}\nolimits}
\def\id{\mathop{\rm id}\nolimits}
\def\Ext{\mathop{\rm Ext}\nolimits}
\def\Ines{\mathop{\rm Ines}\nolimits}
\def\Hom{\mathop{\rm Hom}\nolimits}
\def\bHom{\mathop{\rm\bf Hom}\nolimits}
\def\Rk{ R_\k-\mathop{\bf Mod}}
\def\Rn{ R_n-\mathop{\bf Mod}}
\def\map{\mathop{\rm map}\nolimits}
\def\cf{\mathop{\rm cf}\nolimits}
\def\top{\mathop{\rm top}\nolimits}
\def\Ker{\mathop{\rm Ker}\nolimits}
\def\Bext{\mathop{\rm Bext}\nolimits}
\def\Br{\mathop{\rm Br}\nolimits}
\def\dom{\mathop{\rm Dom}\nolimits}
\def\min{\mathop{\rm min}\nolimits}
\def\im{\mathop{\rm Im}\nolimits}
\def\max{\mathop{\rm max}\nolimits}
\def\rk{\mathop{\rm rk}}
\def\Diam{\diamondsuit}
\def\Z{{\mathbb Z}}
\def\Q{{\mathbb Q}}
\def\N{{\mathbb N}}
\def\bQ{{\bf Q}}
\def\bF{{\bf F}}
\def\bX{{\bf X}}
\def\bY{{\bf Y}}
\def\bHom{{\bf Hom}}
\def\bEnd{{\bf End}}
\def\bS{{\mathbb S}}
\def\AA{{\cal A}}
\def\BB{{\cal B}}
\def\CC{{\cal C}}
\def\DD{{\cal D}}
\def\TT{{\cal T}}
\def\FF{{\cal F}}
\def\GG{{\cal G}}
\def\PP{{\cal P}}
\def\SS{{\cal S}}
\def\XX{{\cal X}}
\def\YY{{\cal Y}}
\def\fS{{\mathfrak S}}
\def\fH{{\mathfrak H}}
\def\fU{{\mathfrak U}}
\def\fW{{\mathfrak W}}
\def\fK{{\mathfrak K}}
\def\PT{{\mathfrak{PT}}}
\def\T{{\mathfrak{T}}}
\def\fX{{\mathfrak X}}
\def\fP{{\mathfrak P}}
\def\X{{\mathfrak X}}
\def\Y{{\mathfrak Y}}
\def\F{{\mathfrak F}}
\def\C{{\mathfrak C}}
\def\B{{\mathfrak B}}
\def\J{{\mathfrak J}}
\def\fN{{\mathfrak N}}
\def\fM{{\mathfrak M}}
\def\Fk{{\F_\k}}
\def\bar{\overline }
\def\Bbar{\bar B}
\def\Cbar{\bar C}
\def\Pbar{\bar P}
\def\etabar{\bar \eta}
\def\Tbar{\bar T}
\def\fbar{\bar f}
\def\nubar{\bar \nu}
\def\rhobar{\bar \rho}
\def\Abar{\bar A}
\def\a{\alpha}
\def\b{\beta}
\def\g{\gamma}
\def\w{\omega}
\def\e{\varepsilon}
\def\o{\omega}
\def\va{\varphi}
\def\k{\kappa}
\def\m{\mu}
\def\n{\nu}
\def\r{\rho}
\def\f{\phi}
\def\hv{\widehat\v}
\def\hF{\widehat F}
\def\v{\varphi}
\def\s{\sigma}
\def\l{\lambda}
\def\lo{\lambda^{\aln}}
\def\d{\delta}
\def\z{\zeta}
\def\th{\theta}
\def\t{\tau}
\def\ale{\aleph_1}
\def\aln{\aleph_0}
\def\Cont{2^{\aln}}
\def\nld{{}^{ n \downarrow }\l}
\def\n+1d{{}^{ n+1 \downarrow }\l}
\def\hsupp#1{[[\,#1\,]]}
\def\size#1{\left|\,#1\,\right|}
\def\Binfhat{\widehat {B_{\infty}}}
\def\Zhat{\widehat \Z}
\def\Mhat{\widehat M}
\def\Rhat{\widehat R}
\def\Phat{\widehat P}
\def\Fhat{\widehat F}
\def\fhat{\widehat f}
\def\Ahat{\widehat A}
\def\Chat{\widehat C}
\def\Ghat{\widehat G}
\def\Bhat{\widehat B}
\def\Btilde{\widetilde B}
\def\Ftilde{\widetilde F}
\def\restr{\mathop{\upharpoonright}}
\def\to{\rightarrow}
\def\arr{\longrightarrow}
\def\LA{\langle}
\def\RA{\rangle}
\newcommand{\norm}[1]{\text{$\parallel\! #1 \!\parallel$}}
\newcommand{\supp}[1]{\text{$\left[ \, #1\, \right]$}}
\def\set#1{\left\{\,#1\,\right\}}
\newcommand{\mb}{\mathbf}
\newcommand{\wt}{\widetilde}
\newcommand{\card}[1]{\mbox{$\left| #1 \right|$}}
\newcommand{\union}{\bigcup}%\limits}
\newcommand{\inters}{\bigcap}%\limits}
\newcommand{\ER}{{\rm E}}
\def\Proof{{\sl Proof.}\quad}
\def\fine{\ \black\vskip.4truecm}
\def\black{\ {\hbox{\vrule width 4pt height 4pt depth
0pt}}}
\def\fine{\ \black\vskip.4truecm}
\long\def\alert#1{\smallskip\line{\hskip\parindent\vrule%
\vbox{\advance\hsize-2\parindent\hrule\smallskip\parindent.4\parindent%
\narrower\noindent#1\smallskip\hrule}\vrule\hfill}\smallskip}

\title{Two Generalizations of co-Hopfian \\ Abelian Groups}
\footnotetext{2010 AMS Subject Classification: Primary 20K10, 20K12; Secondary 20K20, 20K21, 20K30.
Key words and phrases: Abelian groups, directly finite groups, co-Hopfian groups, co-Bassian groups, generalized co-Bassian groups}
\author{Andrey R. Chekhlov \\Faculty of Mathematics and Mechanics, Section of Algebra, \\Tomsk State University, Tomsk, Russia\\{\small e-mails: cheklov@math.tsu.ru, a.r.che@yandex.ru}\\
Peter V. Danchev \\Institute of Mathematics and Informatics, Section of Algebra, \\Bulgarian Academy of Sciences, Sofia, Bulgaria\\{\small e-mails: danchev@math.bas.bg, pvdanchev@yahoo.com}
\\and\\ Patrick W. Keef \\ Department of Mathematics, \\ Whitman College, Walla Walla, WA, Unites States\\{\small e-mail: keef@whitman.edu}}
\maketitle

%\centerline{(To Brendan Goldsmith on his {\bf 75}th birthday)}

\begin{abstract}{By defining the classes of {\it generalized co-Hopfian} and {\it relatively co-Hopfian} groups, respectively, we consider two expanded versions of the generalized co-Bassian groups and of the classical co-Hopfian groups giving a close relationship with them. Concretely, we completely describe generalized co-Hopfian $p$-groups for some prime $p$ obtaining that such a group is either divisible, or it splits into a direct sum of a special bounded group and a special co-Hopfian group. Furthermore, a comprehensive description of a torsion-free generalized co-Hopfian group is obtained. In addition, we fully characterize when a mixed splitting group and, in certain cases, when a genuinely mixed group are generalized co-Hopfian. Finally, complete characterizations of a super hereditarily generalized co-Hopfian group as well as of a hereditarily generalized co-Hopfian group are given, showing in the latter situation that it decomposes as the direct sum of three specific summands.

Moreover, we totally classify relatively co-Hopfian $p$-groups proving the unexpected fact that they are exactly the co-Hopfian ones. About the torsion-free and mixed cases, we show in light of direct decompositions that in certain situations they are satisfactory classifiable -- e.g., the splitting mixed relatively co-Hopfian groups and the relatively co-Hopfian completely decomposable torsion-free groups. Finally, complete classifications of super and hereditarily relatively co-Hopfian groups are established in terms of ranks which rich us that thee two classes curiously do coincide.}
\end{abstract}

\section{Fundamentals and Motivations}

Throughout this article, all groups into consideration are {\it additively} written and {\it Abelian}. Our basic notation and terminology will follow those from \cite{F1,F2}. As usual, for some prime integer $p$, $\Z(p^n)$ denotes the cyclic $p$-group of order $p^n$ for some $n\geq 1$, and $\Z(p^{\infty})$ designates the quasi-cyclic divisible $p$-group.

\medskip

We now proceed with a brief retrospection of our main instruments, which motivate writing up of what is presented in the sequel.

\medskip

A classical notion of some importance in some aspects of abelian group theory is that of a {\it co-Hopfian group} $G$ which means that any injective endomorphism $\phi: G\to G$ is surjective, i.e., $\phi(G)=G$ (equivalently, $\phi$ is an bijective homomorphism, that is, an isomorphism); thus, $\phi(G)$ is a trivial (i.e., a non-proper) direct summand of $G$. In other words, reformulating this in an equivalent manner, $G$ is a group that is not isomorphic to any of its proper subgroups.

\medskip

Some recent non-trivial generalizations of this class of groups were given in \cite{AGS}, \cite{CD1} and \cite{CD2}, respectively.

\medskip

Moreover, recall the classical definition of {\it directly finite groups} that are groups which do not possess an isomorphic proper direct summand. These groups do {\it not} properly retain the property of being co-Hopfian, i.e., there is a directly finite group that is {\it not} co-Hopfian, as well as they are closely related to a generalized variant of Hopfian $p$-groups (see, e.g., \cite{CDGK}).

\medskip

In this vein, we now intend to state the following new concept. It is motivated by the possibility to expand the properties of co-Hopfian groups in light of the aforementioned two classes of groups like this.

\medskip

\noindent{\bf Definition 1.} A group is said to be {\it relatively co-Hopfian} if it is a group that is not isomorphic to a {\bf proper} direct summand of any of its {\bf proper} subgroups.

\medskip

Note that we prefer to use hereafter for further applications the simpler condition stated in Proposition~\ref{altdefRcH}.

\medskip

In fact, such a group cannot be isomorphic to a proper direct summand of any of its proper subgroups, but could be isomorphic to some of its direct summands. So, we may consider two different ways of an investigation, namely of all subgroups (thus including the whole group) and of proper subgroups only (thus excluding the former group).

\medskip

Since the relatively co-Hopfian group does {\it not} have a proper direct summand isomorphic to itself, then such a group does {\it not} possess a non-zero direct summand of the type $A^{(\alpha)}$ for an infinite cardinal $\alpha$. In particular, the divisible groups are relatively co-Hopfian exactly when its torsion-free rank and the rank of each of its $p$-components are all finite; free groups are relatively co-Hopfian exactly when they have finite rank; and in relatively co-Hopfian $p$-groups $G$ every Ulm-Kaplansky invariant $f_n(G)$, defined in the traditional manner as in \cite{F1}, is finite whenever $n<\omega$.

\medskip

On the other hand, mimicking \cite{K}, a group $G$ is called {\it co-Bassian}, provided, for all subgroups $N$ of $G$, whenever $\varphi : G \to G/N$ is an injective homomorphism, then $\varphi(G) = G/N$. In addition, $G$ is {\it generalized co-Bassian} if, for all subgroups $N$ of $G$, whenever $\varphi(G) : G \to G/N$ is an injective homomorphism, then $\varphi(G)$ is a direct summand of $G/N$.

\medskip

These two classes of groups were completely characterized in terms of (generalized) $p$-ranks. Some further extension of the latter class was established in \cite{CDK}.

\medskip

So, motivated by all of this, we now come to the following new concept.

\medskip

\noindent{\bf Definition 2.} We will call a group $G$ {\it generalized co-Hopfian} if, whenever $\phi:G\to G$ is an injective endomorphism, then $\phi(G)$ is a direct summand of $G$.

\medskip

Clearly, invoking to \cite{K}, any co-Bassian group will be co-Hopfian, and a generalized co-Bassian group will be generalized co-Hopfian, but {\it not} necessarily co-Hopfian; for example, any infinite elementary $p$-group is generalized co-Bassian, but will {\it not} be co-Hopfian. Notice also that each divisible group is generalized co-Bassian, but {\it not} co-Hopfian.

\medskip

On the other hand, if we are endowed with the limitation on the direct summand to be {\bf proper} in Definition 2, that is $\phi(G)\not= G$, this version of generalized co-Hopficity somewhat will treat the situation when the group is {\it not} co-Hopfian, which is definitely {\it not} too convenient for our scientific purposes.

\medskip

Moreover, it is also worthy of exploration the more general case when $\phi(G)$ is an essential subgroup in a {\bf proper} direct summand of $G$.

\medskip

Henceforth, our work is organized as follows: In the next Section 2, we present the main results and relevant examples which briefly sound thus: In Propositions~\ref{altdefRcH} and \ref{newcharacter}, we give several characterizations of relative co-Hopficity in certain different aspects. Later on, in Theorem~\ref{torsioncohopfian}, we prove that torsion relatively co-Hopfian groups are exactly the torsion co-Hopfian groups. In Propositions~\ref{alg1} and \ref{compldecomp} we find some characterizing of torsion-free relatively co-Hopfian groups in the cases when they are either algebraically compact or completely decomposable, respectively. In Theorem~\ref{suphered} we determine the structure of super and hereditarily relatively co-Hopfian groups demonstrating that these two classes equal each other. Furthermore, in Proposition~\ref{torfree}, we illustrate that torsion-free generalized co-Hopfian groups are precisely the divisible torsion-free groups, thus exhausting their structure. Next, in Theorem~\ref{ptorsion}, we characterize the structure of generalized co-Hopfian $p$-groups for some arbitrary fixed prime $p$ in terms of direct decompositions of some special divisible groups and co-Hopfian groups. In Proposition~\ref{heregene} and Theorem~\ref{supegene}, we find criteria when a group is hereditarily generalized co-Hopfian and, respectively, super generalized co-Hopfian.

Our further plan to finish the paper is given in Section 3 by raising three significant questions, namely Problems 1, 2 and 3, which solutions hopefully will stimulate a further extremely research study of the current subject.

\section{Principal Results and Examples}

The following general affirmation is very friendly and freely used below.

\begin{lemma}\label{direct} A direct summand of a (relatively, generalized) co-Hopfian group is also so.
\end{lemma}

\Pf Suppose $G=A\oplus B$. Let $\phi:A\to A$ be a monomorphism. Extending $\phi$ by setting it equal to the identity on $B$, it remains injective.

First, if $G$ is relatively co-Hopfian, $\phi(G)=\phi(A)\oplus B$ is essential in $G$ and this immediately implies that $\phi(A)$ is essential in $A$, as needed.

Second, if $G$ is generalized co-Hopfian, it follows that $\phi(G)=\phi(A)\oplus B$ is a direct summand of $G$. If $$G=\phi(G)\oplus C=\phi(A)\oplus B\oplus C,$$ it follow that $$A=\phi(A)\oplus (A\cap (B\oplus C)),$$ as required.
\fine

\subsection{Relatively co-Hopfian Groups}

The following criterion gives a rather more useful characterization of relative co-Hopfian groups than the original treatment, which we will use in the sequel intensively.

\begin{proposition}\label{altdefRcH} The group $G$ is relatively co-Hopfian if, and only if, for every monomorphism $\phi:G\to G$, $\phi(G)$ is essential in $G$. In particular, the group $G\neq\{0\}$ is relatively co-Hopfian if, and only if, all subgroups of $G$ which are isomorphic to $G$ are essential in $G$.
\end{proposition}

\Pf {\bf Necessity.} Let us assume the opposite, namely that there is a monomorphism $\phi:G\to G$ such that $\phi(G)$ is not essential in $G$; that is, there is a non-zero $B\leq G$ such that $\phi(G)\oplus B\leq G$. Since $\phi$ is injective, we obtain
$$\phi(\phi(G))\oplus \phi(B)=\phi(\phi(G)\oplus B)\leq \phi(G)\ne G.$$
As $G\cong \phi(G)\cong \phi(\phi(G))$ and $\phi(B)\ne \{0\}$, it follows that $G$ is isomorphic to a proper direct summand of a proper subgroup of itself, contrary to our initial assumption.

\medskip

{\bf Sufficiency.} It is even more straightforwardly, so we omit the details.

\medskip

The second part follows now immediately.
\fine

The next necessary and sufficient condition provides us with an extra helpful information about the property of a given group to be relatively co-Hopfian. Recall that a group $G$ is known to be {\it directly finite} (or, in other terms, {\it Dedekind finite}), provided $G$ does not possess an isomorphic proper direct summand.

\begin{proposition}\label{newcharacter} A group $G$ is relatively co-Hopfian if, and only if, one of the following three conditions holds:

(i) $G$ is directly finite and the image of any injective endomorphism (= monomorphism) of $G$ is either essential or a proper direct summand.

(ii) There exists a fully invariant essential subgroup of $G$ which is relatively co-Hopfian.

(iii) The inverse image of any non-zero subgroup of $G$ under any injective endomorphism (= monomorphism) of $G$ is non-zero.
\end{proposition}

\Pf (i) Treating the necessity, one elementarily sees that each relatively co-Hopfian group is directly finite; in fact, if we assume the contrary writing $G=A\oplus B$, where $G\cong A$ and $B\neq \{0\}$, then we may write that $G=(A'\oplus B')\oplus B$, where $A=A'\oplus B'$ and $G\cong A'$. So, it must be that $A'\oplus B'\neq G$ and $G\cong A'$ which against our assumption. The other part follows automatically from Proposition~\ref{altdefRcH}.

As for the sufficiency, letting $f:G\to G$ be a monomorphism and $f(G)$ not be essential in $G$. Then, by assumption, one can decomposes $f(G)\oplus A=G$ for some non-zero subgroup $A\leq G$. Hence, there is an obvious isomorphism $G\oplus A  \cong G$ which unambiguously contradicts the direct finiteness of $G$. Thus, $f(G)$ has to be essential in $G$, as required.

(ii) The necessity is very trivial, so we shall be focussed on the sufficiency. Suppose $F$ is a fully invariant essential subgroup of $G$ such that $F$ is relatively co-Hopfian. Letting $f$ be a monomorphism of $G$, one plainly inspects that the restriction $f_F$ is a monomorphism of $F$ and hence $f_F(F)$ is essential in $F$. But, as $F$ is also essential in $G$, we deduce via the transitivity of the essentiality that $f_F(F)$ is too essential in $G$. Consequently, $f(G)$ has to be essential in $G$, as needed.

(iii) Dealing with the necessity, let $M$ be a non-zero subgroup of $G$ and $f$ a monomorphism of $G$. Then, $M\cap f(G)\not=\{0\}$ whence, if $0\not=a\in M$ with $a=f(g)$ for some $g\in G$, then we have that $$0\not=g\in f^{-1}(f(G)\cap M)=f^{-1}(f(G))\cap f^{-1}(M)=G\cap f^{-1}(M)=f^{-1}(M),$$ as we need, where, for any subgroup $N$ of $G$, the set $f^{-1}(N)$ stands for $\{h\in G ~|~ f(h)\in N\}$.

About the sufficiency, if we assume the contrary that there is a monomorphism $f$ of $G$ such that $M\cap f(G)=\{0\}$ for some non-zero subgroup $M$ of $G$, then one follows that $$f^{-1}(f(G)\cap M)=f^{-1}(0)=\{0\},$$ which enables us that $f^{-1}(M)=\{0\}$, thereby contradicting the initial assumption.
\fine

\medskip

This leads us to the question of finding a necessary and sufficient condition when an arbitrary ($p$-)group is directly finite. To this aim, we offer the validity of the following.

\medskip

We remember that a reduced mixed group $G$ with an infinite number of non-zero $p$-primary components $T_p(G)$ is said to be an {\it sp-group}, provided $G$ is a pure subgroup of the Cartesian product $\prod_p T_p(G)$. Notice that the torsion part $T(G)$ of $G$ is the corresponding direct sum $\bigoplus_p T_p(G)$, and the factor-group $G/T(G)$ is always divisible.

\begin{proposition}\label{dirfinsp} The following two statements are true:

(1) A $p$-group $G$ is directly finite if, and only if, the $p$-rank of its divisible part is finite and $f_n(G)$
is finite for all $n<\omega$.

(2) An sp-group $G$ is directly finite if, and only if, its torsion part is directly finite.
\end{proposition}

\Pf The necessity of both (1) and (2) are pretty obvious, so we drop off the details.

We now concentrate on the sufficiency.

(1) Assume that $G=A\oplus C$, where $A\cong G$. Since $$D=D(G)=(D\cap A)\oplus (A\cap C)$$ and $D\cap A\cong D$, it must be that $D=D\cap A$ as $D$ has finite rank meaning that $D\leq A$. Also, since $f_n(G)$ is finite for all
$n<\omega$, some basic subgroup $B$ of the reduced part of $G$ is contained in $A$. So, since $G/B$ is divisible, it follows that $C$ is divisible; but, because the divisible part $D$ of $G$ is contained in $A$, we obtain $C=\{0\}$.

(2) If $G=A\oplus C$, where $A\cong G$, then as above in point (1), $D(G)\leq A$ and $T(G)\leq A$, so $C=\{0\}$, because $G/T(G)$ is divisible.
\fine

We now proceed by including a fast review of some basic ideas regarding essential subgroups. The next technicality is reasonably well-known, so its proof is voluntarily omitted.

\begin{lemma}\label{basic}
Suppose $A\leq G$.

(a) If $G$ is torsion-free, then $A$ is essential in $G$ if, and only if, $G/A$ is torsion if, and only if, $\Q A=\Q G$.

(b) If $G$ is torsion-free of finite rank $n$, then $A$ is essential in $G$ if, and only if, $A$ also has rank $n$.

(c) If $G$ is a $p$-group, then $A$ is essential in $G$ if, and only if, $G[p]\leq A$.

(d) If $G$ is a $p$-group of finite $p$-rank $n$, then $A$ is essential in $G$ if, and only if, $A$ also has $p$-rank $n$.
\end{lemma}

The next statement is somewhat helpful and close to the preceding assertion, so it is worthy of documentation (at least that it seems to be unpublished in the existing literature).

\medskip

Standardly, as noticed above, for the maximal torsion subgroup $T(G)=\oplus_p T_p(G)$ of a group $G$, the symbol $T_p(G)$ stands for its $p$-torsion component for some prime $p$.

\begin{proposition} The subgroup $p^nG$ is essential in $G$ for every integer $n>0$ if, and only if, $T_p(G)$ is divisible. 
\end{proposition}

\Pf {\bf Necessity.} If, in a way of contradiction, we assume that $T_p(G)\neq\{0\}$ is not divisible, then $T_p(G)$ contains a non-zero cyclic direct summand, say $H$, of $G$ of order $k>0$, and hence $H\cap p^{k+1}G=\{0\}$, as pursued.

\medskip

{\bf Sufficiency.} Letting $H$ be a non-zero cyclic subgroup in $G$, if $H\cap T_p(G)\neq\{0\}$, then
$H\cap p^nG\neq\{0\}$ for every integer $n>0$. Furthermore, if $H\cap T_p(G)=\{0\}$, then it is not so difficult to verify that it is possible to choose the subgroup $B$ of $G$ such that $G=T_p(G)\oplus B$ and $H\leq B$ (see, e.g.,
\cite[Theorem 21.2]{F1}). But, since $T_p(B)=\{0\}$, we deduce that $\{0\}\neq p^nH\leq p^nG$ for every $n>0$, as asked for.
\fine

The following is a partial converse to Lemma~\ref{direct}.

\begin{proposition}\label{fullyinv} Let $G=A\oplus B$, where $B$ is fully invariant in $G$, that is, $\Hom(B,A)=\{0\}$. Then, $G$ is relatively co-Hopfian if, and only if, $A$ and $B$ are relatively co-Hopfian.
\end{proposition}

\Pf Certainly, if $G$ is relatively co-Hopfian, then by Lemma~\ref{direct}, both $A$ and $B$ are, as well. So, assume $A$ and $B$ are relatively co-Hopfian, and $\phi:G\to G$ is a monomorphism. It follows that $\phi$ restricts to a monomorphism $B\to B$, which implies that $\phi (B)\subseteq \phi(G)\cap B$ will be essential in $B$. It will, therefore, suffice to show that $\phi(G)\cap A$ is essential in $A$ in order to get applied Proposition~\ref{altdefRcH}.

Certainly, since $\phi(B)\subseteq B$, it induces a homomorphism $\hat \phi: A\cong G/B\to G/B\cong A$. We claim that $\hat \phi$ is also injective: suppose $a\in A$ satisfies $\hat \phi(a)=0$. This means that $\phi(a)\in B$. Since $\phi(B)$ is essential in $B$, there is some $n\in \N$ and $b\in B$ such that $0 \ne n\phi(a)=\phi(b)$; so $b\ne 0$. Since  $\phi(na-b)=0$, and $\phi$ is injective, we can conclude that $na-b=0$, i.e., $0\ne b= na\in A\cap B=0$, which is a contradiction. Therefore, since $A$ is relatively co-Hopfian, we can conclude that $\hat \phi(A)$ is essential in $A$ .

Supposing $0\ne x\in A$, it follows that there is an $n\in \N$ such that $0\ne nx=\hat \phi(a)$, where $a\in A$. In other words, $nx-\phi(a)\in B$. Note that, if $nx=\phi(a)$, we are done. Otherwise, there is an $m\in \N$ and $b\in B$ such that $$mnx-m\phi(a)=\phi(b)\ne 0.$$ Since $b\ne 0$, we can infer that $ma+b\ne 0$, so that $$0\ne \phi(ma+b)=mn x\in \phi(G),$$ as required.
\fine

The following three necessary and sufficient conditions are almost immediate consequences of Proposition~\ref{fullyinv}.

\begin{corollary}\label{divred} Let $G=R\oplus D$ be a group, where $D$ is divisible and $R$ is reduced. Then, $G$ is relatively co-Hopfian if, and only if, $R$ and $D$ are relatively co-Hopfian.
\end{corollary}

\Pf Follows from Proposition~\ref{fullyinv} since $D$ is fully invariant in $G$.
\fine

\begin{corollary}\label{splitmixed} Suppose $G$ is a splitting mixed group, i.e., $G=A\oplus T$, where $T$ is torsion and $A$ is torsion-free. Then, $G$ is relatively co-Hopfian if, and only if, $A$ and $T$ are relatively co-Hopfian.
\end{corollary}

\Pf Follows from Proposition~\ref{fullyinv} since $T$ is fully invariant in $G$.
\fine

\begin{corollary}\label{cotorsionmixed} Suppose $G$ is a reduced co-torsion group with torsion $T$, i.e., $G\cong A\oplus C$, where $A$ is torsion-free and algebraically compact and $C\cong \Ext (\Q/\Z, T)$ is adjusted co-torsion.  Then, $G$ is relatively co-Hopfian if, and only if, $A$ and $C$ are relatively co-Hopfian.
\end{corollary}

\Pf Follows from Proposition~\ref{fullyinv} since $C$ must be fully invariant: to see this, suppose $\phi:C\to A$ is a homomorphism. Since $T$ is torsion and $A$ is torsion-free, we have $\phi(T)=\{0\}$. And since $C/T\cong \Ext(\Q, T)$ is divisible, while $A$ is reduced, we can conclude that $\phi(C)=\{0\}$, as needed.
\fine

Before proceeding by proving our next result on direct sums, we need one more technicality like this.

\begin{lemma}\label{technicalessential} Let $G=A\oplus B$ with projection $\pi: G\to A$, and let $H=K\oplus N\oplus F\leq G$, where $N$ is an essential subgroup of $B$ and the subgroup $\pi(K)$ is essential in $A$. Then, $F=\{0\}$.
\end{lemma}

\Pf Since $N$ is essential in $B$, we have $F\cap B=\{0\}$. We now intend to prove that $F\cap A=\{0\}$, so assume the contrary that $0\neq x\in F\cap A$. Since $\pi(K)$ is essential in $A$, there exists $n\in\mathbb{N}$ such that $nx\neq 0$ and $nx+b\in K$ for some $b\in B$. Note that $b\neq 0$
as for otherwise $0\neq nx\in F\cap K=\{0\}$. If now $0\neq mb\in N$ for $m\in\mathbb{N}$, then $$nmx=(nmx+mb)-mb\in(K\oplus N)\cap F=\{0\}.$$ Therefore, $nmx=0$ and hence $0\neq mb=nmx+mb\in K\cap N=\{0\}$, a contradiction. So, it must be that $F\cap A=\{0\}$.

Assume now that $F\neq\{0\}$ and $0\neq x\in F$. Then, $x=a+b$, where, in view of the already obtained above equalities, $F\cap A=F\cap B=\{0\}$, we have $0\neq a\in A$, $0\neq b\in B$. Thus, there exists $m\in\mathbb{N}$ such that $ma\neq 0$ and $y:=ma+b'\in K$ for some $b'\in B$. Hence, $mx-y=mb-b'\in B$, where $mx-y\neq 0$ since $mx\neq 0$, $y\neq 0$ and $K\cap F=\{0\}$. Consequently, $b_1:=mb-b'\neq 0$, where $b_1\in B$. Observe also that there exists $n\in\mathbb{N}$ such that $0\neq nb_1\in N$. Thus, $$0\neq nb_1=nmx-ny\in (F\oplus K)\cap N=\{0\}.$$ This contradiction gives that $F=\{0\}$, as claimed.
\fine

We now have enough instruments to attack the following which somewhat expands Corollary~\ref{divred} and which gives a new confirmation of the truthfulness of Proposition~\ref{fullyinv}.

\begin{proposition}\label{fullyinv2} Let $G=A\oplus B$, where $B$ is fully invariant in $G$ and both $A$, $B$ are relatively co-Hopfian groups. Then, $G$ also is a relatively co-Hopfian group.
\end{proposition}

\Pf Let $H=K\oplus N\oplus F\leq G$, where $K\cong A$ and $N\cong B$. Since $B$ is fully invariant in $G$, we get $N\leq B$, and since $N\cong B$, we have that $N$ is essential in $B$ thanks to Proposition~\ref{altdefRcH}. So, $K\cap B=\{0\}$ whence, if $\pi: G\to A$ is a projection, then having in mind the  isomorphisms $\pi(K)\cong K\cong A$, one verifies that the subgroup $\pi(K)$ is essential in $A$ by the same Proposition~\ref{altdefRcH}. Consequently, Lemma~\ref{technicalessential} tells us that $F=\{0\}$. So, the group $G$ is really relatively co-Hopfian, as asserted.
\fine

We now proceed by proving a series of some preliminary technicalities.

\begin{lemma}\label{dirsumfull}
Let $G=\bigoplus_{i\in I}~ G_i$, where all $G_i$ are non-zero fully invariant subgroups of $G$. Then, $G$ is a relatively co-Hopfian group if, and only if, every direct component $G_i$ is a relatively co-Hopfian group.
\end{lemma}

\Pf The necessity follows at once applying Lemma~\ref{direct}.

To deal with sufficiency, suppose $\phi:G\to G$ is a monomorphism. It follows that, for all $i\in I$, the map $\phi$ restricts to a monomorphism $G_i\to G_i$. Since each $G_i$ is relatively co-Hopfian, each $\phi(G_i)$ is essential in $G_i$. It follows from general properties of essential subgroups that $\phi(G)=\bigoplus_{i\in I}~ \phi(G_i)$ will be essential in $G=\bigoplus_{i\in I}~ G_i$ too, as asked for in order to apply Proposition~\ref{altdefRcH}.
\fine

As a direct consequence, we derive:

\begin{corollary}\label{torsionreduction}
A torsion group is relatively co-Hopfian if, and only if, each its $p$-component is relatively co-Hopfian.
\end{corollary}

We now state a simple tool for showing a given group fails to be relatively co-Hopfian.

\begin{lemma}\label{infinitesummands} If $G$ is a group with a summand of the form $A:=\bigoplus_{n<\omega} A_n$, where $A_0\ne \{0\}$ and, for each $n<\omega$ there is a monomorphism $\phi_n:A_n\to A_{n+1}$, then $G$ is not relatively co-Hopfian.
\end{lemma}

\Pf Observe elementarily that the maps $\phi_n$s can be joined together to get a monomorphism $\phi:A\to A$ such that $A_0\cap \phi(A)=\{0\}$. This shows that $A$ is not relatively co-Hopfian, so invoking Lemma~\ref{direct} neither is $G$.
\fine

Suppose $A=\oplus_{i\in I} A_i$. It follows from Lemmas~\ref{infinitesummands} and~\ref{basic}(b),(d) that if either (1) each $A_i\cong \Q$, (2) each $A_i\cong \Z$, or (3) there is a prime such that each $A_i\cong \Z(p^\alpha)$, where $\alpha\in \N\cup \{\infty\}$, then $A$ is relatively co-Hopfian if, and only if, $I$ is finite. And using Corollaries~\ref{splitmixed} and \ref{torsionreduction}, this leads to the following:

\begin{corollary}\label{dirsumcyc} Suppose $G$ is either divisible or a direct sum of cyclic groups. Then, $G$ is relatively co-Hopfian if, and only if, it has finite rank and, for all primes $p$, it has finite $p$-rank.
\end{corollary}

Stating now the last assertion in another language, we have:

\begin{corollary}
A direct sum of cyclic groups is relatively co-Hopfian if, and only if, it is Bassian. A divisible group is relatively  co-Hopfian if, and only if, it is co-Bassian.
\end{corollary}

It now follows from Corollary~\ref{dirsumcyc} that there exist direct sums of cyclic groups which are directly finite but {\it not} relatively co-Hopfian, and from Proposition~\ref{dirfinsp} it follows that there exist such sp-groups that are {\it not} relatively co-Hopfian.

\medskip

As usual, for a group $G$, prime $p$ and ordinal $\alpha$, the quotient $U_\alpha^p(G)$ stands for the {\it $p$th Ulm factor} $(p^\alpha G)[p]/(p^{\alpha+1}G)[p]$; we will also let $U_\infty^p(G)=(p^\infty G)[p]$ (see, e.g., \cite{F1,F2}). For $\alpha$ an ordinal or $\infty$, the rank of this quotient is denoted as $f_\alpha^p(G)$ and is called the {\it $\alpha$th Ulm invariant}. We shall omit the prime $p$ from this notation if it is obvious from the context. We also say that $G$ is {\it $p$-semi-standard} if $f_n^p(G)$ is finite for all $n<\omega$.

\medskip

We, thus, arrive at the following.

\begin{corollary}\label{semistandard}
If $G$ is a relatively co-Hopfian group, and $\alpha$ is either $n<\omega$ or $\infty$, then $f_\alpha^p(G)$ is finite for each $p$.
\end{corollary}

\Pf If this failed, $G$ would have a direct summand of the form $\Z(p^{n+1})^{(\omega)}$ or $\Z(p^{\infty})^{(\omega)}$, which cannot happen in view of Lemma~\ref{infinitesummands}, thus substantiating our claim.
\fine

Recall that a subgroup $H$ of a $p$-group $G$ is said to be {\it pure} if, for any $n\geq 1$, the equality $p^nG\cap H=p^nH$ is fulfilled (see, for instance, \cite{F1,F2}). The following is a standard criterion (see, e.g.,  \cite[Section~26(h)]{F1}).

\begin{lemma}\label{purity}
If $G$ is a $p$-group, then $H\leq G$ is pure in $G$ if, and only if, for all $n<\omega$, the equality $(p^n G)[p]\cap H=(p^n H)[p]$ is valid.
\end{lemma}

The following gives a useful tool for $p$-groups that are semi-standard.

\begin{lemma}\label{semistandardendo}
Suppose $G$ is a semi-standard $p$-group. If $\phi:G\to G$ is a monomorphism, then $\phi$ is an isomorphism if, and only if, $\phi(G[p])=G[p]$.
\end{lemma}

\Pf Set $H:= \phi(G)$; so, we need to prove $H=G$. To this end, observe first that since $\phi$ is injective, for all $n\in \N$, we have $\phi((p^n G)[p])=(p^n H)[p]\leq (p^n G)[p]$. In addition, $\phi$ naturally induces a composite of two surjections
\[
G[p]/(p^n G)[p]\to \phi(G[p])/\phi((p^n G)[p])=G[p]/(p^n H)[p]\to G[p]/(p^n G)[p].
\]
However, since $(p^k G)[p]/(p^{k+1} G)[p]$ is finite for all $k<n$, $G[p]/(p^n G)[p]$ is also finite. Thanks to the classical ``Pigeon-hole principle'', we can conclude that this composite is not only a surjection, but a bijection, as well. Therefore, each of the two maps in this sequence must be bijections. Thus, the injectivity if the right morphism ensures that $(p^n G)[p]=(p^n H)[p]$.

\smallskip

Since the latter equality obviously implies $(p^n G)[p]\cap H=(p^n H)[p]$, it follows from Lemma~\ref{purity} that $H$ is a pure subgroup of $G$. However, knowing that $G[p]=H[p]$, \cite[Section~26(j)]{F1} implies that $G=H$, as asked for.
\fine

We are now in a position to prove the following somewhat surprising assertion.

\begin{proposition}\label{pprimary}
If $G$ is a $p$-group, then it is relatively co-Hopfian if, and only if, it is co-Hopfian.
\end{proposition}

\Pf Sufficiency being pretty obvious, suppose $G$ is relatively co-Hopfian. Consulting with Lemma~\ref{semistandard}, $G$ is semi-standard. Letting $\phi: G\to G$ be a monomorphism, it suffices to show $\phi$ is an epimorphism. However, since $G$ is relatively co-Hopfian, $\phi(G)$ must be essential in $G$, so that we must have $G[p]= \phi(G)[p]$, and hence the result follows at once from Lemma~\ref{semistandardendo}.
\fine

In conjunction with this proposition and \cite{BS}, one may ask whether there does exist a co-Hopfian $p$-group that is {\it not} relatively/generalized Hopfian as defined in \cite{CDGK}?

\medskip

Combining Proposition~\ref{pprimary} with Corollary~\ref{torsionreduction}, we immediately deduce that the above result extends to all torsion groups.

\begin{theorem}\label{torsioncohopfian} A torsion group is relatively co-Hopfian if, and only if, it is co-Hopfian.
\end{theorem}

We note the following minor (and well-known) extension of Lemma~\ref{infinitesummands}.

\begin{lemma}\label{torsioncompletesummands} If $G$ is a group with a direct summand $A$ that is an infinite torsion complete $p$-group, then $G$ is not relatively co-Hopfian. In particular, each torsion-complete relatively co-Hopfian $p$-group is finite.
\end{lemma}

\Pf Suppose $A=\overline B$, where $B$ is a basic subgroup of $A$. If $B$ is infinite, then there is a decomposition $B_0\oplus B_1$ with $B_0\ne \{0\}$ and a monomorphism $\phi: B\to B_1$ (as in Lemma~\ref{infinitesummands}, we can think of this as a ``right-shift" operator). This will extend to a monomorphism $$\phi: A=\overline B\to \overline B_0\oplus \overline B_1=A$$ such that $\phi(A)\subseteq \overline B_1$, showing that $A$, and hence $G$, is not relatively co-Hopfian.

The second part is now immediate.
\fine

Recall that a $p$-group $G$ is said to be {\it thick} if whenever $B$ is a direct sum of cyclic $p$-groups and $\phi: G\to B$ is a homomorphism, then $\phi$ is {\it small}, i.e., the kernel of $\phi$ is a {\it large} subgroup of $G$. To use this notion in the present context, we recall the following helpful result.

\medskip

\noindent \cite[Corollary~18(a)]{IK}: {\it If $G$ is a $p$-group that is not thick, then there is a $p$-group $H$ with a direct summand that is an unbounded direct sum of cyclic groups such that $G$ and $H$ embed in each other.}

\medskip

We, thereby, come to the following.

\begin{proposition}\label{coHopfiansarethick} If the $p$-group $G$ is relatively co-Hopfian, then it is thick.
\end{proposition}

\Pf Suppose the contrary that $G$ is not thick and $H=A\oplus B$, where $B$ is an unbounded direct sum of cyclic groups such that $G$ and $H$ embed in each other. Since $B$, and hence $H$, it not relatively co-Hopfian, it follows that $H$ embeds into itself as a non-essential subgroup. Therefore, $G$ also embeds in $H$ as a non-essential subgroup. But since $H$ embeds in $G$, $G$ must embed in itself as a non-essential subgroup. Consequently, $G$ is not relatively co-Hopfian, as stated.
\fine

Note that an unbounded torsion-complete $p$-group will be thick, but not co-Hopfian, so the converse to the last result does not generally hold. It, however, does have the following interesting consequence.

\medskip

Recall that the so-called {\it $\oplus_c$-topology} uses the subgroups $X\leq G$ such that $G/X$ is a direct sum of cyclic groups as a neighborhood base of 0.

\begin{corollary}\label{opluscomplete} Suppose $G$ is a relatively co-Hopfian $p$-group. Then, $G$ is complete in its $\oplus_c$-topology if, and only if, it is finite.
\end{corollary}

\Pf Sufficiency being obvious, suppose $G$ is complete in its $\oplus_c$-topology; in particular, $G$ must be separable. Since $G$ is thick, its $\oplus_c$-completion agrees with its torsion-completion, $\overline G$ (see cf. \cite[Proposition~1.1]{D}). So, Lemma~\ref{torsioncompletesummands} works to derive that $G\cong \overline G$ must be finite, as expected.
\fine

We now mention a specific case of the last above result. In \cite{F2}, L. Fuchs designated the smallest class of Abelian $p$-groups containing the cyclic groups that is closed with respect to direct sums, direct summands, and the torsion subgroups of direct products over non-measurable index sets as the {\it Keef class}, denoting it by $K_p$. Clearly, $K_p$ also contains both the direct sums of cyclic $p$-groups and the torsion complete $p$-groups. It is known that the elements of $K_p$ are all complete in their $\oplus_c$-topologies. Therefore, the following is an immediate consequence of Corollary~\ref{opluscomplete}.

\begin{corollary} A group $G\in K_p$ is relatively co-Hopfian if, and only if, it is finite.
\end{corollary}

We now consider a property that is, in some sense, dual to that of being thick. The $p$-group $G$ is said to be {\it thin} if, for every torsion-complete $p$-group $\overline B$, any homomorphism $\phi: B\to G$ must be small. Ch. Megibben showed that, if $G$ is separable, then $G$ is thin if, and only if, it does not have a subgroup that is isomorphic to an unbounded torsion-complete $p$-group (\cite [Theorem~3.1]{Me}). Thus, the following relates this to our inquiry:

\begin{proposition}\label{coHopfiansarethin}
If $G$ is a separable $p$-group that is relatively co-Hopfian, then $G$ is thin.
\end{proposition}

\Pf We suppose the contrary that $G$ is a separable $p$-group that is not thin, and show that $G$ is also not relatively co-Hopfian. By the aforementioned result, $G$ must have a subgroup $H\cong \overline B$, where $B$ is an unbounded direct sum of cyclic groups. Utilizing now a familiar argument, $\overline B$ has an unbounded proper direct summand, and replacing $H$ with this summand, there is no loss of generality in assuming that $H$ is not essential in $G$.

Note that since the $G$ is semi-standard, if $A$ is a basic subgroup of $G$, then $A$ is also semi-standard, and hence the direct sum of a countable number of cyclic direct summands. By embedding cyclic direct summands of $A$ into larger cyclic direct summands of $B$, there is a monomorphism $\gamma:A\to B$ which extends to a monomorphism $\overline A\to \overline B$. Since $G$ embeds in $\overline A$ which embeds in $\overline B$ which embeds as a non-essential subgroup of $G$ again, we can conclude that $G$ is not relatively co-Hopfian, as stated.
\fine

It follows from a combination of Propositions~\ref{coHopfiansarethick} and \ref{coHopfiansarethin} that a separable relatively co-Hopfian $p$-group is {\it thick-thin}, a class of groups that received considerable attention in \cite{KS}.

\medskip

Recall that the reduced $p$-group $G$ is said to be {\it fully starred} if every subgroup of $G$ has the same cardinality as one of its basic subgroups. For example, if $H$ and $K$ are reduced $p$-groups, it is known that $\Tor (H,K)$ (and all of its subgroups) will be fully starred. In particular, this implies that, for any ordinal $\alpha$, any subgroup of a $p^\alpha$-pure projective $p$-group will be fully starred. In addition, a countable reduced $p$-group is fully starred.

\medskip

We can now record the following.

\begin{proposition} If $G$ is a reduced fully-starred $p$-group, then $G$ is relatively co-Hopfian if, and only if, it is finite.
\end{proposition}

\Pf Certainly, if $G$ is finite, it is co-Hopfian. So, assume $G$ is an infinite reduced and fully-starred co-Hopfian group. Since $G$ must be semi-standard, if it is bounded, it must be finite. Reciprocally, if it is unbounded, then it must have a countably infinite basic subgroup, so that $G$ is also countably infinite and unbounded. Therefore, $G/p^\omega G$ is countably infinite and unbounded, and hence an unbounded direct sum of cyclic groups. But Proposition~\ref{coHopfiansarethick} allows us to infer that $G$, and hence $G/p^\omega G$, must be thick, giving the desired contradiction.
\fine

The following application of the last result on co-Hopficity is well-known (see, e.g., \cite{GG}).

\begin{corollary} A countable reduced $p$-group is relatively co-Hopfian if, and only if, it is finite.
\end{corollary}

We now turn to some results on non-torsion groups. We begin with a simple but applicable observation.

\begin{proposition}\label{TFfactor} Suppose $G$ is a group with torsion $T$. If $T$ and $G/T$ are both relatively co-Hopfian, then $G$ is also relatively co-Hopfian.
\end{proposition}

\Pf Suppose $\phi:G\to G$ is a monomorphism. Clearly, $\phi$ restricts to an injective endomorphism of $T$, and as $T$ is co-Hopfian, as it is relatively co-Hopfian in virtue of Theorem~\ref{torsioncohopfian}, we can conclude that $\phi(T)=T$. This easily guarantees that $\phi$ induces an injective endomorphism $\overline \phi: G/T\to G/T$. Because $G/T$ is relatively co-Hopfian, we can infer that $\overline \phi(G/T)$ is essential in $G/T$. A short argument then forces that $\phi(G)$ is essential in $G$, as required.
\fine

It is worthwhile noticing that, utilizing an analogous idea, the last assertion can be extended thus: {\it If F is a fully invariant subgroup of a group G such that F is co-Hopfian and G/F is relatively co-Hopfian, then G itself is relatively co-Hopfian}.

\medskip

Particularly, we deduce:

\begin{corollary}\label{BcoH}
Suppose $G$ is a group with torsion $T$. If $T$ is co-Hopfian and $G/T$ has finite rank, then $G$ is also relatively co-Hopfian.
\end{corollary}

The last result informs us, for instance, that any Bassian or co-Bassian group is relatively co-Hopfian. The cases of a possible validity of the reverse implication, namely under which extra circumstances $G$ being relatively co-Hopfian yields the same property for either/both $T$ or/and $G/T$, are also of some interest and importance.

\medskip

We now turn to a characterization of the co-torsion-groups that are relatively co-Hopfian. With Corollaries~\ref{divred} and \ref{cotorsionmixed} at hand, we need only consider the cases where $G$ is either (1) reduced, torsion-free and algebraically compact; (2) adjusted co-torsion.

\medskip

We, foremost, start with a useful technical claim.

\begin{proposition}\label{products} If $G=\prod_{i\geq 1}G_i$, where all $G_i$ are non-zero, reduced, torsion-free groups, then $G$ is not a relatively co-Hopfian group.
\end{proposition}

\Pf Let $n_iG_i\neq G_i$ for some integer $n_i>1$. For each $i\geq 1$, let $\phi_i:G_i\to G_i$ be a multiplication by $i\cdot n_i$ and $x\in G_i\setminus n_iG_i$. In the natural way, the $\phi_i$ splice together to an injective homomorphism $\phi:G\to G$; set $H:=\phi(G)$. Considering the vector ${\bf x}:=(x_i)_{i\geq 1}$, it easily follows by a direct check that $H\cap \langle {\bf x}\rangle=\{0\}$, showing that $H$ is not essential in $G$, so that $G$ is not relatively co-Hopfian, as promised.
\fine

Henceforth, for a prime $p$, we denote the ring consisting of all $p$-adic integers by $\mathbb{\widehat{Z}}_p$.

\begin{proposition}\label{alg1} A reduced algebraically compact torsion-free group $G$ is relatively co-Hopfian if, and only if, it is the direct sum of a finite number of copies of $\mathbb{\widehat{Z}}_p$ for various primes $p$.
\end{proposition}

\Pf It is well known that a reduced algebraically compact torsion-free group $G$ is expressible as the direct product $\prod_{p\in \cal P} G_p$, where each component $G_p$ is the $p$-adic completion of a free $\mathbb{\widehat{Z}}_p$-module of $\mathbb{\widehat{Z}}_p$-rank $\kappa_p$, where the $\kappa_p$s are cardinals  (see \cite[Chapter~VII, Proposition~40.1 and Theorem~40.2]{F1}].

Suppose first that $G$ is relatively co-Hopfian. Viewing Corollary~\ref{products}, we can have only a finite number of non-zero terms $G_p$, so that the direct product must be a direct sum.

Note that each term $G_p$ will be fully invariant in $G$, so that we need only consider a particular copy, say $G_p$. We claim that $\kappa_p$ must be finite. To see this, note that if it failed, then $G_p$ would have a direct summand of the form $C$, where $C$ is the $p$-adic completion of $\bigoplus_{n\in \N} \widehat {\Z}_p x_n$. If $\phi:C\to C$ is the monomorphism determines by setting  $\phi(x_n)=p^{n^2}x_n$ for each $n$, then one sees that $C\cong \phi(C)$ and $y=\sum_{n\in \N} p^n x_n$ will satisfy $\phi(C)\cap \langle y\rangle=\{0\}$, showing that $C$ is not relatively co-Hopfian - contradiction.

Conversely, if $G$ is of the form specified, then again, almost all members $G_p$ will be $\{0\}$ and each of them will be fully invariant in $G$. Evidently, if $\phi:G\to G$ is a monomorphism, then $\phi(G_p)\leq G_p$ will be a $\widehat {\Z}_p$-submodule of the same (finite) $\widehat {\Z}_p$-rank as $G_p$. Therefore, the quotient $G_p/\phi(G_p)$ will be a torsion $\widehat {\Z}_p$-module, which means that it is a torsion $p$-group, and thus $\phi(G)$ will be essential in $G$, as it should be.
\fine

We will say the torsion group $T$ is {\it quotient reduced co-Hopfian} if, whenever $\phi:T\to T$ is a monomorphism with $T/\phi(T)$ reduced, then $\phi$ is an epimorphism, i.e., $T/\phi(T)=\{0\}$. It is pretty clear that, if $T$ is co-Hopfian, then it is reduced co-Hopfian. By the usual arguments, a quotient reduced co-Hopfian group cannot have a summand of the form $\Z(p^n)^{(\omega)}$, so that it must be semi-standard. Note also that, if $G$ is a divisible $p$-group of infinite rank, then $G$ will be quotient reduced co-Hopfian, but {\it not} co-Hopfian.

\medskip

The next query arises quite logically.

\medskip

\noindent{\bf Question:} Is it the case that every quotient reduced co-Hopfian torsion group $T$ that is also reduced is co-Hopfian?

\medskip

We can now establish the following.

\begin{proposition}\label{adjusted}
Suppose $G$ is an adjusted co-torsion group; so, if $T$ is the torsion subgroup of $G$, then $G\cong \Ext(\Q/\Z, T)$. Then, the following are equivalent:

(a) $G$ is co-Hopfian;

(b) $G$ is relatively co-Hopfian;

(c) $T$ is quotient reduced co-Hopfian.
\end{proposition}

\Pf It is trivial that (a) implies (b).

To prove the implication (b) implies (c), assume $G$ is relatively co-Hopfian and $\phi:T\to T$ is a monomorphism such that, if $S=\phi(T)$, then $T/S$ is reduced. Observe first that, since $G$ is relatively co-Hopfian, for every prime $p$, the group $G$, and hence the $p$-component of torsion $T_p$, must be $p$-semi-standard. The long-exact sequence for the functors Ext and Hom gives the following:
$$
0=\Hom(\Q/\Z, T/S)\to \Ext(\Q/\Z, S) \to \Ext(\Q/\Z, T),
$$
where $$\Ext (\Q/\Z, S)\cong \Ext(\Q/\Z, T)\cong G.$$ If now $\overline \phi$ is the extension of $\phi$ to $G\to G$, then we can deduce that $\overline \phi$ is injective. It follows that $\overline \phi(G)$ must be essential in $G$. This insures that its torsion subgroup, $S$ say, must be essential in $T$. But, it follows with the aid of Lemma~\ref{semistandardendo} that $S=T$, so that $\phi$ is an epimorphism, as needed.

Assume next that (c) holds and $\overline \phi:G\to G$ is a monomorphism; to establish (a) we need to prove that $\overline \phi$ is surjective. If $\phi$ is $\overline \phi$ restricted to $T$, it follows from the above exact sequence that $T/S$ must be reduced. So, since $T$ is, by assumption, quotient reduced co-Hopfian, we can conclude that $S=T$, i.e., that $\phi$ is an isomorphism. This assures that $\overline \phi$ is also an isomorphism, as required.
\fine

It is worthy of noticing that, if the above Question has an affirmative answer, then the following will be true: let $G$ be a reduced non-zero adjusted cotorsion group. Then, the following three conditions are equivalent: (1) $G$ is a relatively co-Hopfian group; (2) $G$ is a co-Hopfian group; (3) $T=T(G)$ is a co-Hopfian group.

\medskip

%\begin{proposition} The completely decomposable torsion-free group $G$ for which $G=\bigoplus_{t\in\Omega}G_t$, where all direct summands $G_t$ are homogeneous components of $G$ and $\Omega$ is a set of same type direct summands of $G$ of rank~1, is relatively co-Hopfian if, and only if, each $G_t$ has finite rank and each $t\in\Omega$ has only finite type greater than $t$.
%\end{proposition}

%Example: For every prime $p$, let $Q_p\subseteq \Q$ be the fractions of the form $a/p^k$ such that $a\in \Z$, $k\in \N$. Consider $G=\Z\oplus(\bigoplus_{p\in \cal P} Q_p)$. In particular, for all $p\in \cal P$, $Q_p$ has type greater than $\Z$ and there are an infinite number of such $Q_p$.  However, this $G$ pretty clearly is relatively co-Hopfian.  Here is a very slightly different approach:

We now continue with the following description of the relatively co-Hopfian completely decomposable torsion-free groups.

\begin{proposition}\label{compldecomp} Suppose $\{\tau_i\}_{i\in I}$ is a collection of types and, for each $i\in I$, $G_i$ is a non-zero $\tau_i$-homogeneous completely decomposable torsion-free group. Then, $G=\bigoplus_{i\in I} G_i$ is relatively co-Hopfian if, and only if,

(a) Each $G_i$ has finite rank;

(b) There is no infinite subset $\{i_k:k\in \N\}$ such that $\tau_k<\tau_{k+1}$ for each $k\in \N$.
\end{proposition}

\Pf The two conditions being clearly necessary, by Lemma~\ref{infinitesummands}, suppose (a) and (b) both hold; we need to show $G$ is relatively co-Hopfian. Let $\phi:G\to G$ be injective. Note that $\phi$ will extend uniquely to a monomorphism $\Q G\to \Q G$; by Lemma~\ref{basic}(a), it suffices to show that $\phi(\Q G)=\Q G$. If $i\in I$, let $I_i=\{j\in I: \tau_j\geq \tau_i\}$. It now suffices to show that, for all $i\in I$, $$\Q G_i\subseteq \phi(\bigoplus_{j\in I_i} \Q G_j).$$ To this target, we show that, if such an $i\in I$ for which this fails does exist, then there is another index $i'\in I$ such that $i'>i$ and $i'$ also fails to satisfy this condition. Repeatedly applying this idea will produce a counter-example to condition (b), as pursued.

If there is no such $i'>i$, and $J=\{j\in I: \tau_j>\tau_i\}$, then $I_i=\{i\}\cup J$, and $$\bigoplus_{j\in J} \Q G_j= \phi(\bigoplus_{j\in J} \Q G_j).$$ But since $\Q G_i$ has finite dimension and $\phi$ is injective, this ensures that the induced homomorphism
$$\overline \phi: \Q G_i\cong (\bigoplus_{j\in I_j} \Q G_j)/(\bigoplus_{j\in J} \Q G_j)\to (\bigoplus_{j\in I_i} \Q G_j)/\phi(\bigoplus_{j\in J} \Q G_j)\cong \Q G_i$$
is also an isomorphism. This means that
$$\Q G_i\leq \phi(\Q G_i)+ \phi(\bigoplus_{j\in J} \Q G_j)= \phi(\bigoplus_{j\in I_i} \Q G_j),$$
which contradicts the definition of $i$.
\fine

Observe that part (b) of the last result is equivalent to requiring that, for any subset $J\subseteq I$, there is a maximal element $\tau_j\in \{\tau_i:i\in J\}$, i.e., there is not $i\in J$ such that $\tau_i>\tau_j$.

\medskip

We, thus, arrive at the wanted example which shows that the relative co-Hopficity is independent of the ordinary co-Hopficity. Precisely, we show the following.

\medskip

\noindent{\bf Example 1.} There exists a (torsion-free) relatively co-Hopfian group which is {\it not} co-Hopfian.

\medskip

This is pretty straightforward as all co-Hopfian torsion-free groups are themselves divisible, and above, as a source of an example, we exhibited the existence of reduced relatively co-Hopfian torsion-free groups (possibly of infinite torsion-free rank).

Finding a torsion-free relatively co-Hopfian group of finite rank that is not co-Hopfian is also pretty trivial - just let $G$ be any torsion-free group of finite rank that is not divisible.

However, here are other, easier, ways to do this. For instance, just consider a copy of the $p$-adic integers, or consider a torsion-free group whose endomorphism ring is isomorphic to $\Z$; thus, there are such exhibitions of arbitrary rank.

This substantiates our example after all.

\medskip

We are now concerned with the inheritance of the direct sum property of relative co-Hopficity, planning to establish the following statements the first of which generalizes significantly the corresponding result for co-Hopficity.

\begin{proposition}\label{dirsumprop1} If $G$ is a relatively co-Hopfian group and $F$ is a finitely generated group, then $G\oplus F$ too is a relatively co-Hopfian group.
\end{proposition}

\Pf One sees that it is enough to prove the assertion when $F=\langle a\rangle$ is a cyclic group of either infinite order or an order $p^n$ for some prime number $p$ and an integer $n>0$.

To this target, suppose that $A=G\oplus\langle a\rangle$ and $\pi: A\to G$ is the corresponding projection. Assume that $H=C\oplus \langle b\rangle\oplus X\lneqq A$, where $C\cong G$, ${\rm order}(b)={\rm order}(a)$ and $X\neq\{0\}$.

Let ${\rm order}(a)=p^n$. Consider two basic cases as follows:

\medskip

\noindent{\bf Case 1:} Write $b=g+a$ for some $g\in G$ (accurate to certain integer multiple, mutually simple with $p$). We have $A=G\oplus\langle a+g\rangle$ and so $H=(H\cap G)\oplus\langle a+g\rangle$. Here, $H\cap G\cong C\oplus X$ and $H\cap G\lneqq G$ since relatively co-Hopfian groups do not have proper direct summands isomorphic to itself. However, as $C\cong G$ and $X\neq\{0\}$, which is impossible for a relatively co-Hopfian group.

\medskip

\noindent{\bf Case 2:} Write $b=g+p^la$, where $0<l\leq n$. Thus, ${\rm order}(b)={\rm order}(g)$. Since $p^{n-l}b=p^{n-l}g\neq 0$, it must be that $\langle g\rangle\cap(C\oplus X)=\{0\}$. We, thus, arrive at the subgroup $H'=C\oplus\langle g\rangle\oplus X\leq A$, which is a proper direct decomposition.

\medskip

If $C\cap\langle a\rangle=\{0\}$, then $\pi(C)\oplus\pi(\langle g\rangle)\lneqq G$, where $\pi(C)\cong C\cong G$ and $\pi(\langle g\rangle)=\langle g\rangle$, which does not hold for a relatively co-Hopfian group.

Assume now that $C\cap\langle a\rangle\neq\{0\}$, so that $c=p^ta\in C$ for some $0\leq t<n$. Therefore, $C\cap\langle a+g\rangle=\{0\}$. In fact, if $c_1=p^s(a+g)\in C$ for some $0\leq s<n$, then for $m=\max\{t,s\}$ we have $0\neq p^{m-s}c_1-p^ma=p^mg\in\langle g\rangle\cap C$, a contradiction, as expected. Next, if $\pi': A\to G$ is the projection of a group $A=G\oplus\langle a+g\rangle$ on $G$, then
$$\pi'(C\oplus\langle g\rangle)= \pi'(C)\oplus\langle g\rangle\lneqq G,$$ where $\pi'(C)\cong G$, which is manifestly wrong for a relatively co-Hopfian group.

Finally, let ${\rm order}(a)=\infty$. Assume $b=g+ka$ for some integer $k\neq 0$. So, $$H=(H\cap G)\oplus\langle g+ka\rangle\leq A'=G\oplus \langle g+ka\rangle,$$ where $H\cap G\cong C\oplus X$ and $H\cap G\lneqq G$, which is false for a relatively co-Hopfian group.

Let us now we have $b\in G$. If $C\cap\langle a\rangle=\{0\}$, then $\pi(C)\oplus \langle b\rangle\lneqq G$, where $\pi(C)\cong G$, which is untrue for a relatively co-Hopfian group. But, if $C\cap\langle a\rangle\neq\{0\}$, then $c=ka\in C$ for some integer $k\neq 0$. Thus, for $A=G\oplus\langle a+b\rangle$, we obtain $C\cap\langle a+b\rangle=\{0\}$, because if $c_1=k(a+b)$, then $$nc_1-nka=nkb\in \langle b\rangle\cap C=\{0\},$$ a contradiction. So, if $\pi'': A\to G$ is a projection of $A=G\oplus\langle a+b\rangle$ on $G$, then argued as above $\pi''(C)\oplus\langle b\rangle\lneqq G$ and $\pi''(C)\cong G$, which is not fulfilled for a relatively co-Hopfian group.
\fine

As a valuable consequence, we have:

\begin{corollary} If $G$ is a relatively co-Hopfian group and $F$ is a finitely co-generated group, then $G\oplus F$ is too a relatively co-Hopfian group.
\end{corollary}

\Pf Since the finitely co-generated group is known to be a finite direct sum of co-cyclic groups, then the result follows from a combination between Proposition~\ref{dirsumprop1} and Corollary~\ref{divred} equipped with the construction of a divisible relatively co-Hopfian group alluded to above.
\fine

In virtue of \cite{GG1}, a group $G$ is called {\it super relatively co-Hopfian} if every homomorphic image of $G$ is relatively co-Hopfian. Besides, the group $G$ is called {\it hereditarily relatively co-Hopfian} if every subgroup of $G$ is relatively co-Hopfian.

\medskip

We now menage to describe these groups in the next statements.

\begin{lemma}\label{001} Let $G=R\oplus D(G)$, where $D(G)$ is the divisible part of $G$, each $p$-component of $D(G)$ has finite rank, the torsion-free rank of $G$ is finite and, if torsion part of $R$ is non-zero, then each $p$-component is finite. Then, $G$ is super (resp., hereditarily) relatively co-Hopfian.
\end{lemma}

\Pf Note that every subgroup and factor-group of $G$ have the same structure as the group $G$ itself. Suppose $A\leq G$, where $A\cong G$. Thus, $T(G)\leq A$ since the rank of each $T_p(G)$ is finite, and $A$ is essential since the torsion-free rank of $G$ is finite. It now remains to refer to Proposition~\ref{altdefRcH} completing the arguments.
\fine

\begin{proposition} The following are valid:

(1) A divisible group is super (resp., hereditarily) relatively co-Hopfian if, and only if, its torsion-free part and each $p$-component have finite ranks.

(2) A reduced torsion group $G$ is super (resp., hereditarily) relatively co-Hopfian if, and only if, all its
$p$-components are finite.

(3) A reduced torsion-free group $G$ is super (resp., hereditarily) relatively co-Hopfian if, and only if, it has finite rank.

(4) A reduced mixed group $G$ is super (resp., hereditarily) relatively co-Hopfian if, and only if, all its $p$-components are finite and its torsion-free rank is finite.

(5) A non-reduced group $G$ is super (resp., hereditarily) relatively co-Hopfian if, and only if, its divisible and reduced parts are super (resp., hereditarily) relatively co-Hopfian.
\end{proposition}

\Pf In view of Lemma~\ref{001}, we need to prove only the necessity.

(1) It follows from the description of divisible relatively co-Hopfian groups presented above.

(2) Every $p$-component of such a group is super relatively co-Hopfian and, if it is unbounded, then some factor-group is divisible of infinite rank; however, this cannot be happen and so it is bounded. But this means that it has to be finite, because all $f_n(G)$ are finite for $n<\omega$.

For the hereditary case it should be taken into account that any basic subgroup of $G$ is finite as Corollary~\ref{dirsumcyc} claims.

(3) If we assume the contrary that $G$ has infinite rank, then $G$ possesses a reduced factor-group which is unbounded, that certainly cannot be happen for the super case; for the hereditary case $G$ must have a free subgroup of infinite rank, that also cannot be happen, as expected.

(4) It follows from points (2) and (3).

(5) It follows from a combination of (1) to (4) and Proposition~\ref{fullyinv}.
\fine

The next result is subsumed by all of the formulated above being relevant to this matter, arriving at the final assertion of this subsection which extends the corresponding result from \cite{GG1} (compare with \cite{GG} and \cite{GG2} too).

\begin{theorem}\label{suphered} (1) A group $G$ is super relatively co-Hopfian if, and only if, the torsion-free rank and the rank of each $p$-component of divisible part of $G$ are both finite as well as the rank of the reduced part and each $p$-component of the reduced part are also finite.

(2) A group $G$ is hereditarily relatively co-Hopfian if, and only if, it is super relatively co-Hopfian.
\end{theorem}

\subsection{Generalized co-Hopfian Groups}

Recall that, same as in \cite{F1}, the setting $$f_n(G):={\rm rank}(U_n(G))={\rm rank}((p^nG)[p]/(p^{n+1}G)[p])$$ is the {\it $n$th Ulm-Kaplansky invariant} of $G$.

\medskip

Our key statement is the following one. It is elementary, but useful for our further examinations.

\begin{proposition}\label{torfree}
A torsion-free group $G$  is generalized co-Hopfian if, and only if, it is divisible.
\end{proposition}

\Pf It follows straightforwardly by just considering multiplication by all possible integers.
\fine

We are now ready to attack the following criterion in the case of $p$-torsion groups, where $p$ is an arbitrary prime.

\begin{theorem}\label{ptorsion}
A $p$-group $G$ is generalized co-Hopfian if, and only if, either $G$ is divisible, or $G\cong A\oplus C$, where, for some $n<\omega$ and cardinal $\kappa$, $C\cong \Z(p^n)^{(\kappa)}$ and $A$ is a co-Hopfian group such that $f_m(A)=0$ for all $m<n$.
\end{theorem}

\Pf Certainly, if $G$ is divisible, it satisfies either conditions. Otherwise, let $n$ be the smallest positive integer such that $\kappa:= f_{n-1}(G)\ne 0$ and $G=A\oplus C$ be the above decomposition. We need to show such a group $G$ is generalized co-Hopfian precisely when $A$ is co-Hopfian.

Suppose first that $A$ is co-Hopfian and $\phi:G\to G$ is an injection. It follows that $\phi$ restricts to an injection $$p^n G=p^n A\to p^n A=p^n G$$ (see, for instance, \cite{F1}). But, since $A$ is co-Hopfian, it readily follows that $p^n A$ is also co-Hopfian. Therefore, $\phi$ is an automorphism of $p^n G$. It is then plainly argued that $\phi(G)$ is a direct summand of $G$.

\medskip

Conversely, suppose $G$ is generalized co-Hopfian. We first claim that $f_\infty(A)=f_\infty(G)$ is finite, as is $f_{m-1}(A)=f_{m-1}(G)$ whenever $m>n$ is an integer: if, however, this fails, then $G$ has a direct summand of the form $$B:=\Z(p^n)\oplus \Z(p^\alpha)^{(\omega)},$$ where $\alpha$ is either an integer larger than $n$ or $\infty$. There is clearly a monomorphism $$B\to \Z(p^\alpha)^{(\omega)}\leq B.$$ The image of this homomorphism cannot be a direct summand of $B$, since $\Z(p^\alpha)^{(\omega)}$ has no direct summand isomorphic to $\Z(p^n)$. Therefore, since $B$ is not generalized co-Hopfian, owing to Lemma~\ref{direct} neither is $G$.

To show $A$ is co-Hopfian, let $\gamma:A\to A$ be a monomorphism; we now need to establish $\gamma$ is surjective. Exploiting Lemma~\ref{direct}, $A$ will be generalized co-Hopfian, which means that $\gamma(A)$ will be a direct summand of $A$. Consequently, for some subgroup $X\leq A$, we will have $$A=\gamma(A)\oplus X\cong A\oplus X.$$

If, for a moment, $X\ne \{0\}$, then, for some $\alpha\in \omega\cup {\infty}$ with $n<\alpha$, we have $f_\alpha(X)\ne 0$. But this obviously contradicts that $f_\alpha(A)= f_\alpha(A)+f_\alpha(X)$, because all of these values are finite.

Thus, it must be that $X=\{0\}$, which forces that the map $\gamma$ is surjective, as required.
\fine

This assertion simply leads to the following description of the generalized co-Hopfian property in the case of splitting mixed groups.

\begin{corollary}\label{restrdir}
If $G=T\oplus D$ is a group, where $T$ is torsion and $D$ is torsion-free, then $G$ is generalized co-Hopfian if, and only if, $T$ is generalized co-Hopfian (and hence each $p$-primary component $T_p$ of $T$ is generalized co-Hopfian) and $D$ is divisible.
\end{corollary}

\Pf We just apply the second part of Lemma~\ref{direct} in combination with Theorem~\ref{ptorsion} to get the desired claim.
\fine

It was constructed in \cite[Theorem 1.1 (1)(b)]{BS} a Hopfian $p$-group that is {\it not} co-Hopfian. We shall now extend this construction to the following more general setting. Recall that a $p$-group $G$ is said to be {\it semi-standard}, provided all its Ulm-Kaplansky invariants of finite size are finite, that is, for any $n<\omega$, the cardinal number $f_n(G)$ is finite.

\begin{example} There is a Hopfian $p$-group which is {\it not} generalized co-Hopfian. In fact, consider the group mentioned above. We claim that it will {\it not} be generalized co-Hopfian too.

To achieve this, one observes that Theorem~\ref{ptorsion} applies to get that, for semi-standard $p$-groups, co-Hopfian equals to generalized co-Hopfian, and since their group is clearly semi-standard as it is Hopfian, the two notions are obviously equivalent in this context.

This sustains our claim.
\end{example}

Further, one may expect that an analogous claim to that in Proposition~\ref{dirsumprop1} will hold for generalized co-Hopfian groups, namely that if $G$ is a generalized co-Hopfian group and $F$ is a finitely generated group, then $G\oplus F$ is too a generalized co-Hopfian group, but unfortunately this is {\it not} manifestly invalid. In fact, as noticed above, any torsion-free generalized co-Hopfian group is divisible and even if $F$ is taken to be torsion (and hence finite), Theorem~\ref{ptorsion} will enable us that $f_k(F)=0$ for all $k<n$ which is untrue.

\medskip

Any kind of reduced torsion-free groups $G$ for which $\mathrm{E}(G)$ is a subring of the field $\mathbb{Q}$ of all rational numbers clearly contains relatively co-Hopfian groups and also the Hopfian ones. This is an example of relatively co-Hopfian groups that are non co-Hopfian. Such groups can have arbitrary high power.

\medskip

Note that Corollary~\ref{divred} is no longer true for generalized co-Hopfian groups. For example, let $G=D\oplus R$
be such a $p$-group that $D$ is a divisible group of infinite rank and $R\neq\{0\}$ is a generalized co-Hopfian
group of rank $r(R)\leq r(D)$. Then, a routine check shows that $D$ is isomorphic to the divisible hull of $G$, so there exists an injection $G\to D$, but this image is {\it not} a direct summand in $D$, as required.

\medskip

However, the following weaker statement is true:

\begin{lemma}\label{essential} Let $G=A\oplus B$, where the subgroup $B$ is fully invariant in $G$ and, for every homomorphism $f: A\to B$, the $\ker(f)$ is essential in $A$. Then, $G$ is generalized co-Hopfian if, and only if, both $A$ and $B$ are generalized co-Hopfian.
\end{lemma}

\Pf The left-to-right implication follows immediately from Lemma~\ref{direct}.

For the opposite right-to-left implication, suppose that $f: G\to G$ is an injection and $\pi: G\to A$ is the corresponding projection. Since $f(B)\leq B$, it must be that $B=f(B)\oplus B'$. Observe that the restriction $(\pi f)\!\upharpoonright\! A$ is an injection. Indeed, if $0\neq x\in\ker((\pi f)\!\upharpoonright\! A)$, then $f(x)\in B$, but by hypothesis $0\neq mx\in\ker(1_G-\pi)f$ for some integer $m$, so $$0\neq mf(x)=m\pi f(x)\in A\cap B=\{0\},$$ a contradiction.
Hence $A=\pi f(A)\oplus A'$ and, consequently, $$G=\pi f(A)\oplus(f(B)\oplus B'\oplus A').$$

Furthermore, let $\pi_1: G\to \pi f(A)$,$\theta_1: G\to f(B)\oplus B'\oplus A'$ be the corresponding projections.
One sees that the mappings $$\pi f(a)\mapsto f(a)\mapsto (1-\pi)f(a),$$ 
\noindent where $a\in A$ define a homomorphism
$$\varphi:\pi f(A)\to f(B)\oplus B'\leq f(B)\oplus B'\oplus A'.$$

Letting $\pi_1'=\pi_1+\theta_1\varphi_1\pi_1$ and $\theta_1'=\theta_1-\theta_1\varphi\pi_1$, we then calculate that
$$(\pi_1')^2=\pi_1', (\theta_1')^2=\theta_1', \pi_1'+\theta_1'=1_G, \pi_1'\theta_1'=\theta_1'\pi_1'=0_G,$$ whence one may decompose $G=\pi_1'(G)\oplus\theta_1'(G)$. But, since $\im\,\theta_1'\leq\theta_1$ and $$\theta_1'(f(B)\oplus B'\oplus A')=\theta_1(f(B)\oplus B'\oplus A')=f(B)\oplus B'\oplus A',$$ it follows that $\theta_1'(G)=f(B)\oplus B'\oplus A'$.

Next, if $x=\pi f(a)$, where $a\in A$, then $$\pi_1'(x)=x+\theta_1'\varphi(x)=\pi f(a)+(1_G-\pi)f(a)=f(a),$$ and so
$f(A)\leq\pi_1'(G)$. However, since $\pi_1'(G)=\pi_1'(\pi f(A))\leq f(A)$, we have $\pi_1'(G)=f(A)$. Consequently,
$f(G)$ is a direct summand in $G$, as required.
\fine

As for every torsion-free group the kernel of its homomorphisms in a torsion group is already essential, we freely can establish once again the already known Corollary~\ref{restrdir}.

\medskip

An other direct consequence of interest is the following one.

\begin{corollary}\label{fullydir} If $G=\bigoplus_{i\in I}G_i$, where each direct summand $G_i$ is fully invariant in $G$, then $G$ is generalized co-Hopfian if, and only if, each $G_i$ is generalized co-Hopfian.
\end{corollary}

\Pf It follows at once from Lemma~\ref{essential} bearing in mind the fact that the homomorphic image of each $G_i$ is again contained in $G_i$.
\fine

Matching this result with Theorem~\ref{ptorsion}, we may obtain the full determination of the structure of a torsion generalized co-Hopfian group. Moreover, the above corollary has a more transparent verification like this: indeed, after all, if $\phi:G\to G$ is a monomorphism, then apparently each monomorphic image $\phi(G_i)$ will be a direct summand of $G_i$, so that $\phi(G)$ will be a direct summand of $G$ too, as needed.

\medskip

Now, by virtue of \cite{GG1}, a group $G$ is called {\it super generalized co-Hopfian} if every homomorphic image of $G$ is generalized co-Hopfian. Moreover, the group $G$ is called {\it hereditarily generalized co-Hopfian} if every subgroup of $G$ is generalized co-Hopfian.

\medskip

We are now ready to establish the following criterion which expands the corresponding result for hereditary co-Hopficity substantially (see \cite{GG1} for more account).

\begin{proposition}\label{heregene} A group $G$ is hereditarily generalized co-Hopfian if, and only if, it is such a torsion group that $G=A\oplus B\oplus C$, where $C$ is the divisible part of $G$ and each $p$-component of $C$ has finite rank, $B$ is an elementary group, $A$ does not have elementary direct summands and each its $p$-component is finite.
\end{proposition}

\Pf {\bf Necessity.} The torsion property of $G$ follows since the generalized co-Hopfian torsion-free group is always divisible. Furthermore, if either some $p$-component of the reduced part is unbounded or the rank of the divisible
$p$-part is infinite, then $G$ will have a subgroup of the type $\mathbb{Z}(p^n)^{(k)}$ for some $n>1$ and an infinite cardinal $k$. However, this cannot be happen thanking to Theorem~\ref{ptorsion}.

{\bf Sufficiency.} Since each subgroup of $G$ has a structure similar to that of the group itself, Theorem~\ref{ptorsion} applies to get that such a subgroup is generalized co-Hopfian, as required.
\fine

\bigskip

Now, to get characterized the super generalized co-Hopficity, we begin with the following elementary and well-known observation.

\begin{lemma}\label{easy2} Suppose $G$ is a group, $p$ is a prime and $H$ is an epimorphic image of $G$.
If $G$ has finite torsion-free rank and finite $p$-rank, then $H$ has finite $p$-rank.
\end{lemma}

\Pf Setting $L:=G/T$, so $L$ is torsion-free of finite rank. Suppose we are given an epimorphism $G\to H$ and let $K$ be its kernel. Set $K_1=T\cap K$ and $K_2=K/K_1\cong[K+T]/T\leq L$.
It then follows that there is an exact sequence
$$
          K_1/pK_1\to K/pK \to K_2/pK_2.
$$
Since $K_1\leq T$ has finite $p$-rank, so does $K_1/pK_1$. And since $L$ has finite rank, so does $K_2$, and hence $K_2/pK_2$ has finite $p$-rank. Therefore, $K/pK$ also has finite $p$-rank.

Finally, there is another exact sequence
$$
          G[p]\to H[p]\to K/pK,
$$
and since the outside two groups have finite $p$-rank, so does $H$, as stated.
\fine

Unsurprisingly, we also say that the group $G$ is {\it super generalized co-Bassian} if every epimorphic image of $G$ is generalized co-Bassian.

Observe that if $G$ has a summand of the form $\Z(p^k) \oplus \Z(p^\alpha)^{(\omega)}$, where $k\in \N$ and $k<\alpha\in \N\cup \{\infty\}$, it follows from Lemma~\ref{direct} and Theorem~\ref{ptorsion} that $G$ is {\it not} generalized co-Hopfian.

\medskip

We are now attacking the promised above characterization of super generalized co-Hopfian groups along with that of super generalized co-Bassian groups which unexpectedly coincide each other. Specifically, we are establishing the following necessary and sufficient condition which expands the corresponding result from \cite{GG1} (compare with \cite{GG} and \cite{GG2} as well).

\begin{theorem}\label{characterizeSGcB}
Suppose $G$ is a group with torsion subgroup $T$. The following are equivalent:

(a) $G$ is super generalized co-Bassian;

(b) $G$ is super generalized co-Hopfian;

(c) One of following two conditions holds:

\ \ \ (1) $G$ is divisible;

\ \ \ (2) $G/T$ is torsion-free divisible of finite rank and, for each prime $p$, $T_p$ is either divisible, or $pT_p$ has finite $p$-rank.
\end{theorem}

\Pf We start by assuming (a) and proving (b). To that end, if $G$ is super generalized co-Bassian, and $H$ is an epimorphic image of $G$, then $H$ will also be generalized co-Bassian. However, since any generalized co-Bassian group is generalized co-Hopfian, $H$ must also be generalized co-Hopfian. Hence, $G$ is super generalized co-Hopfian.

\medskip

We now assume (b) and prove (c); so assume $G$ is super generalized co-Hopfian. It follows immediately that $G/T$ is torsion-free and generalized co-Hopfian, so it must be divisible. If $T$ is also divisible, then $G$ is divisible and (1) follows. So, assume $p$ is a prime such that $T_p$ fails to be divisible. Thus, $T_p$, and hence $G$, has a direct summand of the form $B\cong \Z(p^k)$ for some $k\in \N$. Write $G=A\oplus B$.

We now claim that $G$ must have finite torsion-free rank: to that target, let $T_A=T\cap A$. We can simply find a free subgroup $F$ of $A$ such that $$A/(T_A\oplus F)\cong (A/T_A)/[F\oplus T_A]/T_A\cong (\Q/\Z)^{(\kappa)}=:X$$ for some infinite cardinal $\kappa$. There is clearly a direct summand $Y$ of $X$ isomorphic to $\Z(p^\infty)^{(\omega)}$, whence it follows that $H:=Y\oplus B$ is an epimorphic image of $G$. However, $H$ is not generalized co-Hopfian, which contradicts that $G$ is super generalized co-Hopfian. Consequently, it follows that $G/T$ is finite-rank divisible.

Let $p$ be any prime such that $T_p\ne 0$ is reduced; we need to show that $pT_p$ has finite $p$-rank. We assume this fails and derive a contradiction, establishing (c).

\medskip

\noindent{\bf Case 1:} The divisible subgroup of $T_p$ has infinite $p$-rank: Fix a decomposition $G=A\oplus B$, where $B\cong \Z(p^k)$ as above. So, $A$ must have a direct summand isomorphic to $\Z(p^\infty)^{(\omega)}$. Therefore, $G$ has a direct summand of the form $\Z(p^k)\oplus \Z(p^\infty)^{(\omega)}$, so that it is not generalized co-Hopfian, contrary to assumption.

\medskip

\noindent{\bf Case 2:} The divisible subgroup of $T_p$ has finite $p$-rank -- note that $T_p/p^2 T_p$ will be a bounded, pure subgroup of $G/p^2 T_p$, and hence it is a direct summand (see \cite{F1} or \cite{F2}). Therefore, $T_p/p^2 T_p$ will be an epimorphic image of $G$. However, if the the divisible subgroup of $T_p$ has finite $p$-rank, but $pT_p$ has infinite $p$-rank, we can conclude that $T_p/p^2 T_p$ must have a direct summand of the form $U:=\Z(p^2)^{(\omega)}$. Thus, $U$ is also an epimorphic image of $G$. This also leads to the fact that $V:=\Z(p)\oplus \Z(p^2)^{(\omega)}$ is an epimorphic image of $G$. But, $V$ is not generalized co-Hopfian, so that $G$ is not super generalized co-Hopfian, contrary to assumption.

\medskip

We now show that (c) implies (a). If (1) holds, then any epimorphic image of $H$ is also divisible, and hence generalized co-Bassian. Suppose, therefore, that (2) holds. Suppose also $\gamma: G\to H$ is an epimorphism and $S$ is the torsion subgroup of $H$; we need to establish that $H$ is generalized co-Bassian. First, one sees that $\gamma$ induces an epimorphism $G/T\to H/S$, so that $H/S$ is also torsion-free divisible of finite rank.

Fix some prime $p$. As $T_p$ is either divisible or $pT_p$ has finite $p$-rank, it follows that $G=A\oplus T_p$ for some $A\leq G$. If $T_p$ is divisible, then $G$ is $p$-divisible. This means that $H$ is $p$-divisible, which, in turn, implies that $S_p$ is $p$-divisible, and hence divisible.

So, assume that $pT_p$ has finite $p$-rank. We claim that $pS_p$ also has finite $p$-rank: note that $\gamma$ reduces to an epimorphism $pG\to pH$. Since $pG$ has both finite torsion-free rank and finite $p$-rank, thanking to Lemma~\ref{easy2} the same holds for $pH$ and, in particular, $pS$ has finite $p$-rank. It, therefore, follows that $S_p$ has, in the language of \cite{K}, generalized finite $p$-rank. So, \cite[Theorem~2.5]{K} applies to deduce that $H$ is generalized co-Bassian, as required.

The proof is now complete after all.
\fine

\begin{corollary}\label{sumssuperGcoB}
Suppose $G=A\oplus B$. Then, $G$ is a super generalized co-Hopfian group (resp., a super generalized co-Bassian group) if, and only if, $A$ and $B$ are both so.
\end{corollary}

The last two statements unambiguously illustrate that the classes of hereditarily generalized co-Hopfian groups and super generalized co-Hopfian groups are different each other, thus contrasting to Theorem~\ref{suphered} quoted above.

\medskip

We are now managing to exhibit certain concrete constructions, thereby demonstrating the independence of the newly defined two group classes. Concretely, we show the following.

\medskip

\noindent{\bf Example 2.} There exists a generalized co-Hopfian $p$-group which is {\it neither} generalized co-Bassian {\it nor} co-Hopfian.

About the first part-half, knowing that there are unbounded co-Hopfian $p$-groups $G$ with finite $f_n(G)$ for all $n\geq 1$ (see, for instance, \cite{BS}), one may write that $|G|=|G/B|$ for some basic subgroup $B$ of $G$ (thus $G/B$ is a divisible quotient). That is why, there is an injection $\eta: G\to G/B$. But one observes that $\eta(G)$ is {\it not} a direct summand in $G/B$, so that $G$ need {\it not} be generalized co-Bassian, as asked for.

The second part-half follows by an immediate application of Theorem~\ref{ptorsion}, substantiating the example.

\medskip

\noindent{\bf Remark.} Nevertheless, with no any difficulty, we may exhibit a reduced mixed group $G$ which is co-Hopfian, but {\it not} generalized co-Bassian.

In fact, let $T$ be an infinite torsion group all $p$-components of which are finite, and set $G:=\overline{T}$ to be the algebraically compact group that is a completion of $T$. Suppose $f: G\to G$ is an injection. Since $f(T)=T$, it must be that $T\leq f(G)$. So, $\overline{T}=G=f(G)$, i.e., $G$ is co-Hopfian.

However, since $T$ is obviously countable, $G$ has the continuum power and the quotient $G/T$ is divisible, there exist two proper subgroups of $G$, say $T<N<G$, such that
$$G/N=D(A)\oplus B, G\cong A ~ {\rm and} ~ B\neq\{0\},$$
where $D(A)$ is the divisible hull of $A$. But $A$ need {\it not} be a proper direct summand of $G/N$, because $G$ is reduced whereas $G/N$ is divisible and $A\neq D(A)$, so that $G$ is really {\it not} a generalized co-Bassian group, as required. The fact that $G$ is not generalized co-Bassian is immediate from \cite[Theorem~2.5]{K}.

\medskip

We now intend to show that, for every prime $p$, each reduced unbounded (and thus infinite) $p$-group of size strictly less than the continuum is neither relatively co-Hopfian nor generalized co-Hopfian. This is closely relevant to \cite[Theorem 1.1 (2)]{BS}. Indeed, about relative co-Hopficity, we may directly apply Proposition~\ref{pprimary} to get the desired result. As for the generalized co-Hopficity, Theorem~\ref{ptorsion} works thus: if we assume the contrary that there is such a reduced unbounded $p$-group $G$ of size strictly less than the continuum that is generalized co-Hopfian, then $G$ has to be co-Hopfian with finite Ulm-Kaplansky invariants in some segment. This, however, contradicts the fact that was established in \cite{BS}.

Note that in Theorem~\ref{ptorsion} the direct summand $C$ can have an arbitrary cardinality, and the other direct summand $A$ could be finite and even zero, so that the requirement "unbounded" is essential and cannot be replaced by the more weak choice of being "infinite".

So, the Remark now sustained.

\section{Concluding Discussion}

Here we are now concerned with some additional queries as the major one is whether the classes of generalized co-Hopfian and relatively co-Hopfian groups are independent each other. The answer is definitely {\bf yes} as the next things show:

\medskip

$\bullet$ ~ One knows by what we have shown above that all generalized co-Hopfian torsion-free groups are divisible, but there are non-divisible relatively co-Hopfian groups, which illustrates that the converse cannot happen.

\medskip

$\bullet$ ~ There is a generalized co-Hopfian $p$-group which is {\it not} relatively co-Hopfian -- indeed, we just need to take the cardinal $\kappa$ to be infinite in Theorem~\ref{ptorsion}. Moreover, any divisible group of infinite rank is generalized co-Hopfian but {\it not} relatively co-Hopfian.

On the other side, any reduced generalized co-Hopfian $p$-group $G$ having finite Ulm-Kaplansky invariants $f_n(G)$ for all $n<\omega$ has to be co-Hopfian and so relatively co-Hopfian -- indeed, we may just invoke Theorem~\ref{ptorsion}.

\medskip

We end the work with the next three intriguing problems. They are closely connected to the procedure described in \cite[Section 4]{BS}, where it is demonstrated that if a $p$-group $G$ is co-Hopfian, then so does its first Ulm factor $G/p^{\omega}G$ under the presence of the Martin's axiom. We thus arrive at the following two challenging questions, the first one of which is pertained to Proposition~\ref{dirfinsp}(1) and the second one to Theorem~\ref{ptorsion}, but both of them are results in ZFC without any extra set-theoretic assumptions.

\medskip

\noindent{\bf Problem 1.} If a $p$-group $G$ is directly finite, is its first Ulm factor $G/p^{\omega}G$ also directly finite?

\medskip

\noindent{\bf Solution.} The question follows directly from Proposition~\ref{dirfinsp}, because it must be that $f_n(G/p^{\omega}G)=f_n(G)$ for $n<\omega$ and $f_{\alpha}(G/p^{\omega}G)=0$ for $\alpha\geq\omega$ (see \cite{F1}). This ends our proof.

\medskip

Note that, if a $p$-group $G$ is relatively co-Hopfian, then its first Ulm factor $G/p^{\omega}G$ is also relatively co-Hopfian. In fact, just a simple combination of the conclusions in \cite[Section 4]{BS} and Proposition~\ref{pprimary} works to get the claim.

\medskip

\noindent{\bf Problem 2.} If a $p$-group $G$ is generalized co-Hopfian, is its first Ulm factor $G/p^{\omega}G$ also generalized co-Hopfian?

\medskip

\noindent{\bf Solution.} Thankfully, Theorem~\ref{ptorsion} yields that $G=A\oplus C$, where $C$ is a $p^n$-bounded group and $A$ is a co-Hopfian group. So, one writes that $$G/p^{\omega}G\cong (A/p^{\omega}A)\oplus C.$$

Furthermore, since by what we have commented above $A/p^{\omega}A$ is co-Hopfian, and $f_m(A/p^{\omega}A)=f_m(A)$
for every $m<\omega$ (see \cite{F1}), Theorem~\ref{ptorsion} is workable to get that $G/p^{\omega}$ is generalized co-Hopfian, as pursued. The proof is over.

\medskip

So, we come to our final question which settling is unknown to us yet.

\medskip

\noindent{\bf Problem 3.} Find suitable conditions on the first Ulm subgroup $p^{\omega} G$ of a given $p$-group $G$ under which, if the first Ulm factor $G/p^{\omega}G$ is directly finite (resp., generalized co-Hopfian), then so does the whole group $G$.

\medskip

\noindent{\bf Funding:} The work of the first-named author, A.R. Chekhlov, is supported by the Ministry of Science and Higher Education of Russia (agreement No. 075-02-2024-1437).

\end{document}